\documentclass[12pt]{amsart}
 \usepackage{
 amsmath,amsfonts,
 amssymb,amsthm
}
 \numberwithin{equation}{section}
 \newtheorem{prop}{Proposition}[section]

 \newtheorem{cor}{Corollary}[section]

 \newtheorem{thm}{Theorem}[section]
  \newtheorem*{thmm}{Theorem}
  \newtheorem*{conj}{The Hessian Conjecture}

 \newtheorem{dfn}{Definition}[section]

 \newcommand{\p}{\partial}

 \newcommand{\sA}{\mathcal A} 
 
 \newcommand{\sG}{\mathcal G} 
  \newcommand{\sI}{\mathcal I}
 \newcommand{\sC}{\mathcal C}

 \newcommand{\sL}{\mathcal L}

 \newcommand{\sP}{\mathcal P}

 \newcommand{\al}{\alpha}
 
 \newcommand{\de}{\delta}
 \newcommand{\ep}{\varepsilon}

 \newcommand{\De}{\Delta}
 \newcommand{\Ga}{\Gamma}

 \newcommand{\C}{\mathbb C}

 \newcommand{\tr}{\operatorname{tr}}

\newcommand{\Aut}{\operatorname{Aut}}

\newcommand{\sBB}{\mathcal B}
\newcommand{\TTT}{\mathbb T}

\begin{document}
\title{Colored graphs, Gaussian integrals and stable graph polynomials.}
  \author{I.V.Artamkin}
  \address{State University Higher School of Economics, Moscow}
  \thanks{Supported by the grants of RFBR   08-01-00110
of RFBR and SU HSE 09-01-12185-ofi-m (09-09-0010) and 09-01-0063,
and by the grant of the Laboratory of mathematical investigations
TZ-62.0(2010) and by the  project   "Development of new methods of
study of integrable systems and moduli spaces in geometry topology
and mathematical physics" of the Federal Special Program.}
  \email{artamkin@mail.ru}
\date{}

 \maketitle



Asymptotic expansions of  Gaussian integrals may often be
interpreted as generating functions for certain combinatorial
objects (graphs with additional data). In this article we discuss a
general approach to all such cases using colored graphs. We prove
that the generating power series for such graphs satisfy the same
system of partial differential equations as the Gaussian integral
and the formal power series solution of this system is  unique.
 The solution is obtained as the   genus
expansion  of the generating power series. The initial term of this
expansion is the corresponding generating function for trees. The
consequence   equations for this term turns to be equivalent to the
inversion problem for the gradient mapping defined by the initial
condition. The equations for the higher terms of the genus expansion
are linear. The solutions of these equations can be expressed
explicitly by substitution of the initial conditions and the initial
term (the tree expansion) into some universal polynomials (for
$g>1$) which are generating functions for   stable closed graphs.
(For $g=1$ instead of
 polynomials appears logarithm.) The stable graph polynomials satisfy
 certain recurrence. In \cite{A} some of these results were obtained for $r=1$ by
 more or less direct solution of  differential equations. Here we present
 purely combinatorial proofs.

\section{Introduction.}

First let us fix the notations. Consider vectors $X=(x_1,  \ldots ,
x_r)$ and a symmetric $r\times r$ matrix $S=(s_{ij})$. We shall
consider  $x_i$ and $s_ {ij}$ and $\hbar $ as  independent
commutative variables.

 For a multi-index $N=(n_1,  \ldots ,n_r)  $
 we shall use
the notations $X^N=x_1^{n_1} \cdot \ldots \cdot x_r^{n_r}   $,
$N!=n_1!  \ldots n_r!  $ and $|N|=\sum n_i$. $N\ge 0$ will mean that
all $n_i\ge 0$. For the multi-index $ (0, \ldots ,0,1,0, \ldots ,0)$
(all zeros except 1 in the position $i$) we shall use the notation
$\{ i \}$:
\begin{equation}\label{def_i}
  \{ i \} =(0, \ldots ,0,1,0, \ldots ,0).
\end{equation}
Next we denote $
  \{ i j \} =   \{ i  \} +  \{ j \}$,
   so that for $i\ne j$ the multi-index $\{ i j\}$ has exactly two
non-zero positions $i$ and $j$ and for $i= j$ the multi-index $\{ i
i\}$ has $2$ in the non-zero positions $i$. Thus for $i\ne j$  $\{ i
j\} != 1$ and $\{ i i\} !=2$. As well we may use multi-indices
 $ \{ i j  k\} =   \{ i  \} +  \{ j \} +  \{ k \}$ and so on.

We start with a formal power series $U(X,\hbar )\in \frac 1{\hbar
}\C [[ X, \hbar ]]$ which we shall  write
 as a Taylor expansion
\begin{equation}\label{U_expan}
  U(X, \hbar )=
  \sum _{g \ge 0 } \sum _{N \ge 0}
     a_{g,N}
     \frac {X^N}{N!}  \hbar ^{g-1}.
     \end{equation}
It is convenient to consider the coefficients $ a_{g,N}$ as
independent variables.

A formal definition of modular graph see in \cite{A}. Informally
speaking modular graph is a graph which may have edges with only one
end.  Such edges are usually called {\it tails}. Now we want to
label vertices of a modular graph by by the variables $ a_{g,N}$,
tails will be labeled by $x_1, \ldots , x_r$ and the edges will be
labeled by the elements of the symmetric matrix $S$. For that
purpose we need to mark the half-edges incident to each vertex by
the numbers $1, \ldots ,r$ so that two half-edges of one edge may be
marked by different indices. Thus it is natural to insert a new
two-valent vertex into the middle of the edge which will break the
edge into two new edges. Now each of them can be marked by one index
$i \in \{ 1, \ldots ,r \}$. This leads to  the following formal
definition.

\begin{dfn} \label{Def_bip_graph} A
bipartite colored modular  graph is a collection of the following
data:
\begin{enumerate}
  \item a modular graph $\Ga $ whose set of vertices  $V(\Ga )$
   is a
  disjoint union of two partite sets $V(\Ga )=V_a(\Ga )\sqcup
V_s(\Ga )  $;
 \item a mapping from the set of edges and tails of $\Ga $ to the
 set $  \{ 1, \ldots ,r \}$ which we shall call coloring;
 \item a mapping $g:V_a(\Ga )\to \{ 0,1,2,\ldots \}$; nonnegative
 integer g(v) will be called   genus of the vertex $g$. We shall call
 a graph combinatorial if for all its
 vertices $g(v)=0$.
   \end{enumerate}
These data should satisfy the following properties:
\begin{enumerate}
     \item two vertices from the same partite sets are not connected by an edge;
    \item vertices from $ V_s(\Ga )  $ should be only two-valent and should have
   no incident tails.
   \end{enumerate}
The  vertices from $V_s(\Ga )$ will be called $s$-vertices and the
 vertices from $V_a(\Ga )$ will be called $a$-vertices.
\end{dfn}

 Consider a bipartite colored modular  graph $\Ga $. {\it Genus} of
  $\Ga $ is defined by
\begin{equation}\label{def_genus}
g(\Ga )=\sum _{v\in V_a(\Ga )} g(v) + b_1(\Ga )- b_0(\Ga )+1,
\end{equation}
where  $b_m(\Ga )$ is
 the $m$-th Betti number of the graph (considered as a $1$-dimensional
 simplicial complex). Thus for a connected graph $\Ga \ $
\begin{equation}\label{def_genus_conn}
g(\Ga )=\sum _{v\in V_a(\Ga )} g(v) + b_1(\Ga )).
\end{equation}
For a graph $\Ga $ let us fix a multi-index
   $N(\Ga )=(n_1,  \ldots ,n_r) $
  where   $n_i$  are   the numbers of    and tails
   colored by $i$.  A graph without tails (i.e.  $N(\Ga )=(0,  \ldots ,0) $)   will
 be called {\it closed}. {\it Valence} of  an $a$-vertex $v$ of a colored graph $\Ga $
is a multi-index $N(v)=(\nu _1(v), \ldots  ,     \nu _r(v) ) $,
where
  $\nu _i(v)$ is
 the number of   edges and tails colored by $i$ and incident to $v$.

Now to each $a$-vertex  $v\in V_a(\Ga )$ we attach the variable
 $  a_{g(v),N(v)}   $ and we attach $s_{ij} $ to an $s$-vertex  $v\in V_s(\Ga )$
if $v$ is incident to two edges colored by $i$ and $j$. Thus for
each colored graph $\Ga $ we define the monomial
\begin{equation}\label{mu(Ga)}
 \mu(\Ga )= \prod _{v\in V_a(\Ga )} a _{g(v),\nu (v)}  \prod _{v\in V_a(\Ga )} s_{ij}.
\end{equation}
For $N=(n_1,  \ldots ,n_r) $, $g\ge 0$  denote by $\tilde \sBB
_{g,N} $ the set of all colored bipartite graphs $\Ga $ having $n_i$
 half-edges of the color $i$ and $g(\Ga )=g$ and  by $ \sBB _{g,N} $ the set of
all such connected bipartite  graphs. Put $\tilde \sBB _{g}= \bigcup
_{N\ge 0} \tilde \sBB _{g,N} $ and $ \sBB _{g}= \bigcup _{N\ge 0}
\sBB _{g,N} $;  $\tilde \sBB   = \bigcup _{g\ge 0} \tilde \sBB _{g}
  $ and   $ \sBB   = \bigcup _{g\ge 0}   \sBB _{g}   $.

Define two generating series
\begin{equation}\label{Psi_def_int}
    \Psi (S,X, \hbar )=\sum _{\Ga \in  \sBB }
\frac 1{|\Aut \Ga |}
    \mu (\Ga )  X^{N(\Ga )}  \hbar ^{g(\Ga )-1} \in   \frac 1{\hbar }\C
[[ X, S, \hbar ]] ;
\end{equation}
\begin{equation}\label{tilde_Psi_def_int}
   \tilde  \Psi (S,X, \hbar )=\sum _{\Ga \in   \tilde  \sBB }
  \frac 1{|\Aut \Ga |}
    \mu (\Ga )  X^{N(\Ga )}  \hbar ^{g(\Ga )-1} \in   \frac 1{\hbar }\C
[[ X, S, \hbar , \frac 1{\hbar }]] .
\end{equation}
A standard combinatorial principle says that the generating function
for all graphs is the exponent of the generating function for
connected graphs.

\begin{thm}
\begin{equation}\label{exp_Psi}
\tilde \Psi (X,S, \hbar )=\exp\left[ { \Psi (S,X, \hbar )} \right] .
\end{equation}
\end{thm}
Note that substituting $S=0$ into (\ref{Psi_def_int}) we get the
generating power series for distinct $a$-vertices i.e.
\begin{equation}\label{Psi_S=0}
  \Psi ( X,0, \hbar )=
  \sum _{g \ge 0 } \sum _{N \ge 0}
     a_{g,N}
     \frac {X^N}{N!}  \hbar ^{g-1}=  U(X, \hbar ).
\end{equation}

We prove that the generating series $\tilde \Psi $ is the unique
solution of the system of linear partial differential equations
generalizing the heat equation and the series $  \Psi $ is the
unique solution of the system of nonlinear partial differential
equations generalizing the Burgers equation (see \cite{A}).

\begin{thm}\label{th_main_1}
1)  The generating series (\ref{tilde_Psi_def_int}) for all colored
modular graphs
$$   \tilde  \Psi (S,X, \hbar )=\sum _{\Ga \in   \tilde  \sBB }
  \frac 1{|\Aut \Ga |}
    \mu (\Ga )  X^{N(\Ga )}  \hbar ^{g(\Ga )-1} \in   \frac 1{\hbar }\C
[[ X, S, \hbar , \frac 1{\hbar }]]  $$
 is the unique solution in $\C[[ x, S, \hbar
,
 \frac 1{ \hbar } ]]$ of the equations
\begin{equation}\label{Main_sys_Burg_Psi_tilde}
\frac {\p \tilde  \Psi }{\p s_{ij}}=
 \frac {\hbar }{\{ ij \}  !}    \frac {\p ^2 \tilde \Psi }{\p x_i \p x_j}
   \end{equation}
 with the initial condition
 \begin{equation}\label{Init_cond_Burg_Psi_tilde}
   \tilde   \Psi (X, 0, \hbar )= \exp \left[ U(X,  \hbar )
  \right]
\end{equation}
 and $\tilde \Psi $
provides the formal asymptotic expansion of the Gaussian integral:
\begin{equation}\label{asymp_exp}
\tilde\Psi ( x, S, \hbar ) \sim \frac 1{ (2\pi \hbar)^{r/2} (\det S)^{1/2}}\int \exp
\left[
   U(\xi ,\hbar) -\frac 1{ 2 \hbar } (X-\xi )^T S^{-1}  (X-\xi )\right] d
   \xi
\end{equation}
   \\
   2) The generating series (\ref{Psi_def_int}) for all colored
modular graphs
$$\Psi (S,X, \hbar )=\sum _{\Ga \in  \sBB } \frac 1{|\Aut \Ga |}
    \mu (\Ga )  X^{N(\Ga )}  \hbar ^{g(\Ga )-1} \in   \frac 1{\hbar }\C
[[ X, S, \hbar ]]$$
 is the unique solution in $\frac 1{ \hbar } \C[[ x, S, \hbar  ]]$ of the
 equations
 \begin{equation}\label{Main_sys_Burg_Psi}
  \frac {\p  \Psi }{\p s_{ij}}=
 \frac {\hbar }{\{ ij \}  !}   \left[ \frac {\p ^2 \Psi }{\p x_i \p x_j}
+ \left( \frac {\p  \Psi }{\p x_i} \right )
 \left( \frac {\p  \Psi }{\p x_j} \right )    \right]
 \end{equation}
 with the initial condition
    \begin{equation}\label{Init_cond_Burg_Psi}
    \Psi (X, 0, \hbar )=  U(X,  \hbar )
 \end{equation}
 and $\Psi $
provides the formal asymptotic expansion of the Gaussian integral
 \begin{equation}\label{asymp_exp_log}
  \Psi ( X, S, \hbar ) \sim \log  \frac 1{ (2\pi \hbar)^{r/2} (\det S)^{1/2}}\int \exp
\left[
   U(\xi ,\hbar) -\frac 1{ 2 \hbar } (X-\xi )^T S^{-1}  (X-\xi )\right] d
   \xi
   \end{equation}
\end{thm}

Note that the statements concerning Gaussian integrals are more or
less trivial. For a positive matrix $S$ the integral
(\ref{asymp_exp}) is simply the average of the function $U$ by the
normal distribution having mean value $X$ and covariance matrix
$\hbar S$. It is easy to verify that this integral considered as a
function on $X$ and $S$ satisfies  the system
(\ref{Main_sys_Burg_Psi_tilde}). Thus its asymptotic expansion
should also satisfy it. Now it is clear that the formal power series
solution of (\ref{Main_sys_Burg_Psi_tilde}) and
(\ref{Main_sys_Burg_Psi}) with the corresponding initial condition
exist and is unique.

Next let us consider the genus expansion of $\Psi $
\begin{equation}\label{gen_exp_Psi_bip}
\Psi (S,X,\hbar )= \sum _{g\ge 0} \Psi _g (S,X)\hbar ^{g-1},
\end{equation}
where $ \Psi _g (S,X) \in \C [[X,S]]$.

 For $g=0$   (\ref{Main_sys_Burg_Psi}) provides the equation
\begin{equation}\label{Burg_Psi_0_bip}
\frac {\p  \Psi _0}{\p s_{ij}}=
 \frac {1 }{\{ ij \}  !}   \left( \frac {\p  \Psi _0}{\p x_i} \right )
 \left( \frac {\p  \Psi  _0}{\p x_j} \right )    .
\end{equation}
For $g>0$   (\ref{Main_sys_Burg_Psi}) provides recursive equations
\begin{equation}\label{Burg_Psi_g_bip}
\frac {\p  \Psi _g}{\p s_{ij}}=
 \frac 1{\{ ij \}  !}   \left[ \frac {\p ^2 \Psi  _{g-1}}{\p x_i \p x_j}
+\sum _{m=0}^g \left( \frac {\p  \Psi _m}{\p x_i} \right )
 \left( \frac {\p  \Psi   _{g-m}}{\p x_j} \right )    \right] .
\end{equation}
The  equation  (\ref{Burg_Psi_0_bip}) looks better for the gradient
vector function
\begin{equation}\label{def_Phi_bip}
\Phi (S,X)= \nabla   _X \ \Psi _0(S,X )= \left( \Phi _1 (S,X),
\ldots ,\Phi _r(S,X)\right) .
\end{equation}
where $\Phi _i(S,X)= \frac {\p  \Psi _0(S,X)}{\p x_i} $. Then
(\ref{Burg_Psi_0_bip}) provides the equations:
\begin{equation}\label{Burg_Phi_bip}
\frac {\p  \Phi _m}{\p s_{ij}}=  \Phi _i \frac {\p  \Phi _m}{\p x_j}
+\Phi _j   \frac {\p  \Phi }{\p x_i} \quad \mbox{for} \quad   i\ne
j;
\end{equation}
\begin{equation}\label{Burg_Phi_bip=}
\frac {\p  \Phi _m}{\p s_{ii}}= \Phi _i \frac {\p  \Phi _m}{\p x_i}
\quad \mbox{for}  \quad   i= j.
\end{equation}

 Initial conditions for these equations are
given by the genus expansion of the initial condition
(\ref{Init_cond_Burg_Psi}):
\begin{equation}\label{U_expans}
  U(X, \hbar )=
  \sum_{g\ge 0}U_g(X) \hbar ^{g-1}.
     \end{equation}
We shall also use the gradient vector function
\begin{equation}\label{F_def}
  F(X)=
   \nabla  _X \  U_0(X) = \left( F_1(X), \ldots ,  F_r(X) \right) ,
     \end{equation}
where  $F_i  = \frac {\p    U _0}{\p x_i } $ and the Hessian matrix
function
\begin{equation}\label{H_def}
  H(X)=
   \nabla  _X \  F(X) =
  \left( \frac {\p ^2  U _0}{\p x_i\p x_j} \right) .
     \end{equation}
The equations (\ref{Burg_Phi_bip})--(\ref{Burg_Phi_bip=}) may be
solved explicitly. They provide the following functional equation.
\begin{thm} \label{th_fun_eq}
Let $F(X) = \left( F _1(X), \ldots ,F _r(X)\right)$ be any formal
series vector.  Then the solution of the system (\ref{Burg_Phi_bip})
and  (\ref{Burg_Phi_bip=}) with the initial condition $\Phi (
0,X)=F(X)$ satisfies
 the functional equation
\begin{equation}\label{fun_eq_Phi_bip}
\Phi(S,X)=F(X+S\Phi(S,X))  .
\end{equation}
\end{thm}
It is not hard to verify this by direct calculations but we give in
section \ref{sect_genus exp} a combinatorial proof for this theorem.

The functional equation (\ref{fun_eq_Phi_bip}) is equivalent to the
inversion problem for the formal mapping $A(X)=  X- SF(X)$. This
well-known fact was discussed for the case of diagonal matrix $S$ in
\cite{Asia}. Here we shall only present the statement; the proof is
quite the same as in \cite{Asia}.

\begin{cor}\label{inv_AB}
Let $\Phi (S,X)$  be the solution of  the system
(\ref{Burg_Phi_bip}) and  (\ref{Burg_Phi_bip=}) with the initial
condition $\Phi ( 0,X)=F(X)$. Consider the following
 formal mappings from $\C ^r$ to  $\C ^r$:
  \begin{equation}\label{inv1_Phi}
  A(X)=  X- SF(X)
\end{equation}
\begin{equation}\label{inv2_Phi}
  B(X)= X+ S\Phi(S,X)
\end{equation}
Then these
  mappings are inverse to each other:
 \begin{equation}\label{inv}
    A(B(X))=X \mbox {    and     }\  B(A(X))=X.
\end{equation}
\end{cor}

Thus the formal Cauchy problems for the systems of the Burgers
equations for $\Phi $ is equivalent to the problem of finding the
inverse function for the initial conditions. Integrating  $\Phi $ we
may get the first term of the expansion  (\ref{gen_exp_Psi_bip}). It
is remarkable that the second term of these
 expansions may be presented explicitly in terms of $\Phi $ and $H$.

 The  equation  (\ref{Burg_Psi_g_bip})   provide the following
 equations for $\Psi _1 $:
\begin{equation}\label{Burg_Psi_1_bip}
\frac {\p  \Psi _1}{\p s_{ij}}=
 \frac 1{\{ ij \}  !}   \left[ \frac {\p ^2 \Psi  _{0}}{\p x_i \p x_j}
+        \Phi _i
   \frac {\p  \Psi   _{1}}{\p x_j}
+        \Phi _j
   \frac {\p  \Psi   _{1}}{\p x_i}
      \right] .
\end{equation}

\begin{thm}\label{th_Phi_1}
The solution of the  equations (\ref{Burg_Psi_1_bip}) or
 with the initial conditions $ \Psi _1
(0,X)=U_1(X)$  is given by:
\begin{equation}\label{Phi_1}
 \Psi _1 (S,X) = U_1 \left( X+ S\Phi(S,X)\right) -
 \frac 12 \tr \ln \left( E- S H(X+ S\Phi(S,X))\right)
 \end{equation}
 \end{thm}
Probably it is still possible to prove this theorem by direct
calculation but we present a combinatorial proof of it in section
\ref{sect_genus exp}.

For $g>1$ the     recurrent  equations (\ref{Burg_Psi_g_bip})   for
$\Psi _g $ are linear on  $\Psi _g $. The solution may be expressed
 in terms of $\Phi $, $H$ and the {\it stable graph polynomials}.
A graph is called {\it stable} if all its genus zero vertices are at
least trivalent and    all its genus one vertices are at least
univalent. In  section  \ref{sect_genus exp}  we introduce for $g>1$
{\it stable graph polynomials} $P_g (\{ a_{g,N}\} , S)$ (see
(\ref{stable_Psi})) depending on independent variables $ a_{g,N} $
for all $g\ge 0 \ $  $|N|\ge 3$ and   symmetric matrix $ S$ as
generating functions for stable graphs. The stable graph polynomials
satisfy certain recurrence (see theorem \ref{prop_rec_P_g}) and
certain homogeneity properties (see theorem
\ref{homogen_st_gr_poly}). The solution of (\ref{Burg_Psi_g_bip}) is
expressed by the stable
 graph polynomials as follows (see theorem \ref{prop_Psi_via_P_g}):
 \begin{multline}\label{Psi_via_P_g_1}
    \Psi _g ( S,X )= \\
=P_g\left( \left\{  a_{g,N}:= \frac {\p ^{|N|}  U _g \left(X + S
\Phi ( S,X ) \right) }
     {\p X ^N } \right\} ,\left ( E-  SH \left(   X+
 S \Phi ( S,X )\right) \right) ^{-1}  S \right) ;
\end{multline}

We may arrange the stable graph polynomials into the generating
power series
\begin{equation}\label{stab_pol_gen}
   \sP(\{ a_{g,N}\} , S , \hbar  ) = \sum _{g\ge 2} P_g (\{ a_{g,N}\} , S) \hbar ^{g-1}.
\end{equation}
Substituting $X=0$ into (\ref{asymp_exp_log}) we obtain another
useful asymptotic expansion.
\begin{thm}\label{th_asymp_exp_stab}
If $a_{g,N}=0$ for $|N|+2g-2\le 0$ then the seriees
(\ref{stab_pol_gen}) provides the asymptotic expansion
\begin{equation}\label{asymp_exp_stab}
  \sP (\{ a_{g,N}\} , S , \hbar  ) \sim \log  \frac 1{ (2\pi \hbar)^{r/2} (\det S)^{1/2}}\int \exp
\left[
   U(\xi ,\hbar) -\frac 1{ 2 \hbar }  \xi  ^T S^{-1}   \xi  \right] d
   \xi ,
   \end{equation}
   where  $U(X, \hbar )=
  \sum _{g \ge 0 } \sum _{N \ge 0}
     a_{g,N}
     \frac {X^N}{N!}  \hbar ^{g-1}$.
\end{thm}

Note that the case $r=1$ which is far from being trivial. For this
case we have one variable $s$ corresponding to edges of a graph and
two-index variables $a_{g,n}$. Denoting by $\sA _g^k$ the set of
genus $g$ stable closed  graphs  we define
\begin{equation}\label{stable_P_g_r=1_intro}
 P _g(\{  a_{m,N} \} ,s)=\sum _{k=0}^{3g-3}
 \sum _{\Ga  \in  \sA _g^k } \frac {\mu (\Ga )}{|\Aut \Ga |}s^k .
\end{equation}
(Stable genus $g$ graph without half-edges has at most $3g-3$
edges.) For instance for $g=2$
\begin{multline}\label{P_2_1_intro}
  P_2  = a_{2,0}+ \frac 12
  a_{1,1 } ^2  s +  \frac 12
  a_{1,2 }  s
  +\frac 12
  a_{1,1  } a_{0,3  } s^2+\frac 18 a_{0,4  } s^2+
 \frac 5{24}
  a_{0,3  } ^2 s^3.
\end{multline}

There are two interesting specializations of the variables $\{
a_{g,n} \}$: counting functions for all combinatorial graphs and
counting functions for all stable combinatorial graphs and

For the counting functions for all combinatorial graphs
  we put
\begin{equation}\label{all_gr_mu_intro}
    a _{g,n}^{\mathrm{comb}} = \left\{\begin{array}{cc}
      1  & \mbox{ if  }g=0  \\
      0 & \mbox{ otherwise, }  \\
     \end{array}\right. \end{equation}
and for the counting functions for all stable combinatorial graphs
  we put
\begin{equation}\label{all_gr_mu_intro}
    a _{g,n}^{\mathrm{st}} = \left\{\begin{array}{ccc}
      1  & \mbox{ if  }g=0 &  \mbox{ and  } n\ge 3  \\
      0 & \mbox{ otherwise. } & \\
     \end{array}\right. \end{equation}
The function  $\Phi $ satisfies the functional equation
\begin{equation}\label{all_gr_fun_eq_intro}
\Phi ^{\mathrm{comb}}(s,x)=e^{ x+ s \Phi ^{\mathrm{comb}}(s,x)}
 \end{equation}
for the counting functions for all combinatorial
graphs\footnote{Note that $ \Phi ^{\mathrm{comb}}(s,0)$ is the
classical generating function for rooted trees (without half-edges)
whose coefficients are given by the well-known Caley formula:
$$  \Phi ^{\mathrm{comb}}(s,0)=\sum _{k=0}^{\infty} \frac{(k+1)^k}{(k+1)!}s^k. $$
 } and
 the functional equation
 \begin{equation}\label{all_gr_st_fun_eq_intro}
e^{ x+ s \Phi ^{\mathrm{st}}(s,x)}=1+ x+( s+1) \Phi
^{\mathrm{st}}(s,x)
 \end{equation}
for the counting functions for all stable combinatorial graphs.

Stable graph polynomials for both cases coincide; we denote
\begin{equation}\label{P_g^comb_intro}
   P_g^{\mathrm{comb}}(s)= P _g\left( \left\{  a_{m,N}:=
    a_{m,N}^{\mathrm{comb}}   \right\}  \right)  =
  P _g\left( \left\{  a_{m,N}: =  a_{m,N}^{\mathrm{st} }  \right\} ,
 s \right) .
 \end{equation}
For instance for $g=2$
\begin{equation}\label{P_2_comb_intro}
  P_2 ^{\mathrm{comb}} = \frac 18  s^2+
 \frac 5{24}
  s^3.
\end{equation}

In section \ref{r=1} we prove the following theorem.
\begin{thm} \label{count_comb_st}
1) Combinatorial stable graph polynomials  $P_g^{\mathrm{comb}}$ for
$g>2$ satisfy the recurrence
\begin{multline}\label{rec_P_g_r=1_count_intro}
    \frac {d P_g^{\mathrm{comb}}}{d s}= \frac 12 \left[
  D_{comb}^2( P _{g-1}^{\mathrm{comb}})
   +2s D_{comb}( P _{g-1}^{\mathrm{comb}}) +
 \sum _{m=2}^{g-2} D_{comb}( P^{\mathrm{comb}} _{m})D_{comb}( P^{\mathrm{comb}} _{g-m})
 \right] ,
\end{multline}
where
\begin{equation}\label{D_comb_intro}
   D_{comb}= s(s+1)\frac{d}{ds}  - (g-1).
\end{equation}
2) For $g\ge 2$ the counting function for all combinatorial graphs
\begin{equation}\label{stable_P_g_r=1_intro}
 \Psi _g ^{\mathrm{comb}}(s,x )
= \frac 1{\Phi ^{\mathrm{comb}}(s,x )^{g-1}}
 P _g ^{\mathrm{comb}}\left(
 \frac {s\Phi ^{\mathrm{comb}}(s,x )}{1- s \Phi ^{\mathrm{comb}}(s,x ) }\right).
\end{equation}
3) For $g\ge 2$ the counting function for all stable combinatorial
graphs
\begin{equation}\label{stable_P_g_r=1_st_intro}
 \Psi _g ^{\mathrm{st}}(s,x )=
 \frac 1{\left( 1+x+(s+1) \Phi ^{\mathrm{st}}(s,x) \right)^{g-1}}
  P _g  ^{\mathrm{comb}}  \left(
 \frac {s\left( 1+x+(s+1) \Phi ^{\mathrm{st}}(s,x) \right)}
 {1-s (x+(s+1)\Phi ^{\mathrm{st}}(s,x ) ) }\right)
 .
\end{equation}
\end{thm}
Formula (\ref{stable_P_g_r=1_st_intro}) was derived in \cite{A} by
direct solution of the equations (\ref{Burg_Psi_g_bip}).

\section{Bipartite colored graphs.} \label{bip_graph}

First let us prove that $\tilde \Psi (X,S, \hbar )$ satisfies the
system (\ref{Main_sys_Burg_Psi_tilde}). For this purpose we need to
interpret the derivatives of  $\tilde \Psi (X,S, \hbar )$ and
$\tilde \Psi (X,S, \hbar )$ as certain generating power series.

Define the sets of bipartite colored modular (connected) graphs with
one marked $ij$-valent $s$-vertex by $\tilde \sBB _{g,N}^{[ij]}$ ($
\sBB _{g,N}^{[ij]}$). For a graph $\Ga \in \tilde \sBB
_{g,N}^{[ij]}$ (or $ \sBB _{g,N}^{[ij]}$) definition of $\mu (\Ga )$
should be improved: we put
\begin{equation}\label{mu(Ga)_bip_mar}
 \mu^{[ij]}(\Ga )= \prod _{\begin{array}{c}
    \mbox{\footnotesize nonmarked} \\
  {
  v\in V_s(\Ga ) }\\
  \end{array}
   } s _{ij}
  \prod _{v\in V_a(\Ga )} a _{g(v),\nu (v)}.
\end{equation}
Fix $L=(l_1, \ldots , l_r)$,  $N=(n_1, \ldots , n_r)$. The set of
bipartite colored modular (connected) graphs with $n_i+l_i$ tails of
the color $i$, $l_i$ of them marked and ordered, will be denoted by
$\tilde \sBB _{g,N,[L]}$ ($ \sBB _{g,N,[L]}$).

It is clear that
\begin{equation}\label{marked_ed_Psi_tild_bip}
\sum _{-\infty <g<+\infty  } \sum _{N \ge 0}
    \left( \sum _{\Ga \in \tilde  \sBB _{g,N}^{[ij]} }
    \frac {\mu (\Ga )}{|\Aut \Ga |}\right)
     X^N \hbar ^{g-1}=
     \frac {\p  \tilde \Psi (S,X,\hbar )}{\p s_{ij}  },
\end{equation}
\begin{equation}\label{marked_half_Psi_tild_bip}
\sum _{-\infty <g<+\infty  } \sum _{N \ge 0}
     \left( \sum _{\Ga \in \tilde \sBB _{g,N,[L]}^{K} }
    \frac {\mu (\Ga )}{|\Aut \Ga |}\right)
     X^N \hbar ^{g-1}=
     \frac {\p ^l \tilde \Psi  (S,X,\hbar )}{\p x_1^{l_1} \ldots \p x_r^{l_r} },
\end{equation}
and the same is true for the generating series for  connected graphs
(i.e. for $\Psi $  without tilde):
\begin{equation}\label{marked_ed_Psi_bip}
\sum _{g\ge 0  } \sum _{N \ge 0}
    \left( \sum _{\Ga \in   \sBB _{g,N}^{[ij]} }
    \frac {\mu (\Ga )}{|\Aut \Ga |}\right)
     X^N \hbar ^{g-1}=
     \frac {\p   \Psi (S,X,\hbar )}{\p s_{ij}  },
\end{equation}
\begin{equation}\label{marked_half_Psi_bip}
\sum _{g\ge 0  } \sum _{N \ge 0}
     \left( \sum _{\Ga \in  \sBB _{g,N,[L]} }
    \frac {\mu (\Ga )}{|\Aut \Ga |}\right)
     X^N \hbar ^{g-1}=
     \frac {\p ^l  \Psi  (S,X,\hbar )}{\p x_1^{l_1} \ldots \p x_r^{l_r} },
\end{equation}

There is a natural clutching  map $ \tilde \sBB _{g,N,[\{ ij \} ]}
\to \tilde \sBB _{g+1,N}^{[ij]}$: to clutch together two tails we
insert a $ij$-valent $s$-vertex between them. This map is bijective
for $i\ne j$ and is $2$-fold covering for  $i= j$ (since there are
two possible orderings on the set of two marked tails). Hence the
generating power series coincide up to factor $\frac {\hbar }{\{ ij
\} !}$. So we get exactly the equation
(\ref{Main_sys_Burg_Psi_tilde}). Thus we have proved the  theorem
\ref{th_main_1}.

\section{Genus expansion.} \label{sect_genus exp}

Next let us consider the genus expansion of $\Psi $
(\ref{gen_exp_Psi_bip}):
$$ 
\Psi (S,X,\hbar )= \sum _{g\ge 0} \Psi _g (S,X)\hbar ^{g-1},
$$ 

As we have seen in the Introduction
 for $g=0$  the equations (\ref{Main_sys_Burg_Psi}) provide  the
 equations (\ref{Burg_Psi_0_bip})
$$ 
\frac {\p  \Psi _0}{\p s_{ij}}=
 \frac {1 }{\{ ij \}  !}   \left( \frac {\p  \Psi _0}{\p x_i} \right )
 \left( \frac {\p  \Psi  _0}{\p x_j} \right )    .
$$ 
For $g>0$   (\ref{Main_sys_Burg_Psi}) provides recursive equations
(\ref{Burg_Psi_g_bip})
$$ 
\frac {\p  \Psi _g}{\p s_{ij}}=
 \frac 1{\{ ij \}  !}   \left[ \frac {\p ^2 \Psi  _{g-1}}{\p x_i \p x_j}
+\sum _{m=0}^g \left( \frac {\p  \Psi _m}{\p x_i} \right )
 \left( \frac {\p  \Psi   _{g-m}}{\p x_j} \right )    \right] .
$$ 
For the gradient vector functions
$$ 
\Phi (S,X)= \nabla   _X \ \Psi _0(S,X )= \left( \Phi _1 (S,X),
\ldots ,\Phi _r(S,X)\right) ,
$$ 
$$ 
  F(X)=
   \nabla  _X \  U_0(X) = \left( F_1(X), \ldots ,  F_r(X) \right) ,
   $$ 
 and the Hessian matrix
function
$$ 
  H(X)=
   \nabla  _X \  F(X) =
  \left( \frac {\p ^2  U _0}{\p x_i\p x_j} \right) .
 $$ 
we got the equations (\ref{Burg_Phi_bip}) and (\ref{Burg_Phi_bip=})
$$ 
\frac {\p  \Phi _m}{\p s_{ij}}=  \Phi _i \frac {\p  \Phi _m}{\p x_j}
+\Phi _j   \frac {\p  \Phi }{\p x_i} \quad \mbox{for} \quad   i\ne
j;
$$ 
$$ 
\frac {\p  \Phi _m}{\p s_{ii}}= \Phi _i \frac {\p  \Phi _m}{\p x_i}
\quad \mbox{for}  \quad   i= j.
$$ 

Let us prove the theorem  \ref{th_fun_eq}:
\begin{thmm} {\bf 1.3}
Let $F(X) = \left( F _1(X), \ldots ,F _r(X)\right)$ be any formal
series vector.  Then the solution of the system (\ref{Burg_Phi_bip})
and  (\ref{Burg_Phi_bip=}) with the initial condition $\Phi (
0,X)=F(X)$ satisfies
 the functional equation
\begin{equation}\label{fun_eq_Phi_bip}
\Phi(S,X)=F(X+S\Phi(S,X))  .
\end{equation}
\end{thmm}

 Note that according to (\ref{marked_half_Psi_bip})
the series
 $\Phi _i (S ,X)$ is the generating function for the genus $0$ connected  trees
from $ \sBB _{0, [\{ i\} ]}= \bigcup _{N\ge 0} \sBB _{0,N, [\{ i \}
]}$ (trees with one marked tail of the color $i$). Consider the set
$\sBB _{0, [\{ i\} ]}^0$ of edgeless genus $0$ connected trees with
one marked half-edge of the color $i$ $\sBB _{0, [\{ i\} ]}^0$ (i.
e. $\sBB _{0, [\{ i\} ]}^0$ is the set of single vertices).
Attaching to each graph $\Ga \in  \sBB _{0, [\{ i\} ]}$ the vertex
adjacent to the marked tail provides the mapping
\begin{equation}\label{map_c0_bip_trees}
   c_0 \colon  \sBB _{0, [\{ i\} ]} \to  \sBB _{0, [\{ i\} ]}^0.
\end{equation}
For any $\De \in  \sBB _{0, [\{ i\} ]}^{(0)}$ ($\De $ consists of a
single vertex with several tails and one marked tail of the color
$i$)
 all the graphs in  $c_1^{-1} (\De )$
are constructed from $\De $ by clutching arbitrary genus $0$ trees
with one marked tail to some of the tails of $\De $ (inserting a
two-valent $s$-vertex  between the corresponding tail of $\Ga $ and
the marked edge of the tree). Therefore
\begin{equation}\label{c0^(-1)_vertex}
          \sum _{\Ga  \in   c_0^{-1} (\De )  }
    \frac {\mu (\Ga )}{|\Aut \Ga |}
     X^{N(\Ga )}
=     \frac {\mu (\De )}{|\Aut \De |}
   \left(    X +  S \Phi (S,X )\right) ^{N(\De )}.
\end{equation}
But the generating function for $ \sBB _{0, [\{ i\} ]}^{(0)}$ is
$\Phi _i (0,X)=F_i(X)$ and therefore taking the sum over all $\De
\in  \sBB _{0, [\{ i\} ]}^{(0)}$ we get
\begin{multline}\label{c0^(-1)_vertex1}
  \Phi _i (S,X ) = \sum _{\Ga  \in  \sBB _{0, [\{ i\} ]}   }
    \frac {\mu (\Ga )}{|\Aut \Ga |}
     X^{N(\Ga )}=\sum _{\De \in \sBB _{0, [\{ i\} ]}^{(0)}}
          \sum _{\Ga  \in   c_0^{-1} (\De )  }
    \frac {\mu (\Ga )}{|\Aut \Ga |}
     X^{N(\Ga )} = \\
=  \sum _{\De \in \sBB _{0, [\{ i\} ]}^{(0)}}   \frac {\mu (\De
)}{|\Aut \De |}
   \left(    X +  S \Phi (S,X )\right) ^{N(\De )}
   =F_i   \left(    X +  S \Phi (S,X )\right) .
\end{multline}

The functional equation (\ref{fun_eq_Phi_bip}) is equivalent to the
inversion problem for the formal mapping $A(X)=  X- SF(X)$.

\begin{cor}\label{inv_AB}
Let $\Phi (S,X)$  be the solution of  the system
(\ref{Burg_Phi_bip}) and  (\ref{Burg_Phi_bip=}) with the initial
condition $\Phi ( 0,X)=F(X)$. Consider the following
 formal mappings from $\C ^r$ to  $\C ^r$:
  \begin{equation}\label{inv1_Phi}
  A(X)=  X- SF(X)
\end{equation}
\begin{equation}\label{inv2_Phi}
  B(X)= X+ S\Phi(S,X)
\end{equation}
Then these
  mappings are inverse to each other:
 \begin{equation}\label{inv}
    A(B(X))=X \mbox {    and     }\  B(A(X))=X.
\end{equation}
\end{cor}

Differentials of inverse mappings are also inverse to each other.
Consider the Hessian matrix
\begin{equation}\label{def_Hess}
\Theta  (S,X)= \nabla   _X \ \Phi (S,X)=
 \left( \frac {\p ^2 \Psi _0 }{\p x_i\p x_j} \right)  .
\end{equation}
Note that   $\Theta (0,X)$  is the Hessian matrix
of the initial condition 
 (\ref{H_def}):     $\Theta (0,X) =H(X)$. Thus we get the following equation for
$\Theta $.
\begin{cor}
\begin{multline}\label{inv_Teta}
E+ S \Theta(S,X) = \left( E- S \Theta\left( 0,X+ S\Phi(S,X)\right)
\right) ^{-1}= \\ = \left( E- S H\left( X+ S\Phi(S,X)\right) \right)
^{-1}
\end{multline}
\end{cor}

The following considerations will provide us an independent
combinatorial proof of corollary \ref{inv_Teta}.

An $s$-vertex $v'\in V_s(\Ga )$ of a colored bipartite connected
modular graph $\Ga \in  \sBB _{g,N}$ will be called
 $1$-{\it cut} if deletion  of $v'$
 disconnects the graph and at least one of the two new connected components
 has genus zero. A graph without $1$-cuts will be called $2$-{\it connected}.
 Pick a graph $\Ga \in  \sBB _{g,N}$.
 Let us mark all the vertices  $v''\in V_a(\Ga )$ which are  connected by an edge
 with
 at least one
vertex $v'\in V_s(\Ga )$ which is not a $1$-cut
 and all the the vertices $v''\in V_a(\Ga )$
 with $g(v'')>0$. (Note that a genus $0$ graph will have no marked vertices.)
 Next let us delete all the $0$-cuts  connected by an edge
 with
 at least one  marked vertex. As the
   result we shall obtain a number of
 genus $0$ connected component $\Ga _i$, $i>0$  and
  one genus $g$ connected component $\Ga _1$
   without  $1$-cuts (the one having marked vertices).
Let us for $g>0$ denote the set of all colored bipartite connected
and $2$-connected (i.e.  without $1$-cuts) modular graphs $\Ga \in
\sBB _{g,N}$  by $ \sBB _{g,N}^{(1)}$.
   The above construction provides the mapping
\begin{equation}\label{map_c0_bip}
   c_1 \colon  \sBB _{g,N} \to  \sBB _{g,N}^{(1)}.
\end{equation}
Consider for  $g>0$ the corresponding generating function
\begin{equation}\label{Psi_def_bip_1}
    \Psi _g ^{(1)}(S,X )= \sum _{N \ge 0}
     \left( \sum _{\Ga \in \sBB _{g,N}  ^{(1)}}
    \frac {\mu (\Ga )}{|\Aut \Ga |}\right)
     X^N
\end{equation}
For any $\Ga \in  \sBB _{g,N}^{(1)}$ all the graphs in  $c_1^{-1}
(\Ga )$ may be constructed from $\Ga $ by clutching arbitrary genus
$0$ trees with one marked tail to some of the tails of $\Ga $
(inserting a two-valent $s$-vertex  between the corresponding tail
of $\Ga $ and the marked edge of the tree). Therefore
\begin{equation}\label{c1^(-1)_Ga}
          \sum _{\De  \in   c_1^{-1} (\Ga )  }
    \frac {\mu (\De )}{|\Aut \De |}
     X^{N(\De )}
=     \frac {\mu (\Ga )}{|\Aut \Ga |}
   \left(    X +  S \Phi (S,X )\right) ^N.
\end{equation}
Summing over all genus $g$ $2$-connected graphs we obtain the
following expression of $\Psi _g $ via  $\Psi _g ^{(1)}$
\begin{prop}\label{prop_Psi_via_Psi_1} For $g>0$
\begin{equation}\label{Psi_via_Psi_1}
    \Psi _g (S,X )=  \Psi _g ^{(1)}(S,X+ S \Phi (S,X )) .
\end{equation}
\end{prop}

    Consider the set $ \sL _{N,[ij]}^{k} $ of colored bipartite connected modular
trees $\Ga $ consisting of a chain of $k$ genus $0$ vertices $V_a
(\Ga )$
 interleaving
with $k-1$ two-valent $s$-vertices $v'\in V_s(\Ga )$ having two
marked ordered tails $i$ and $j$ incident to the farthest vertices
of \  $\Ga $. Put $\sL _{N,[ij]}=\bigcup _{k\ge 0} \sL _{N,[ij]}^{k}
$ and $\sL _{[ij]}=\bigcup _{N\ge 0} \sL _{N,[ij]} $ . Note that
such graphs have no nontrivial automorphisms (at least for the case
of ordered tails);  reading all the vertices of $\Ga $ along the
chain starting from the  vertex incident to the first marked tail
matches $\mu (\Ga )\frac {X^N}{N!}$ to a certain summand of the $ij$
element of the matrix
\begin{equation}\label{chain}
    \underbrace {HS H S H \ldots
      S H  }_{k \mbox{  times  } H},
\end{equation}
where $H=H(X)$
   is the   the Hessian matrix $(\frac {\p ^2 U _0( X)}{\p x_i  \p x_j })$.
Consider the matrix $\Upsilon _k(X)$ of generating series defined by
\begin{equation}\label{Ups_def}
 \Upsilon  _k(X)_{ij}=\sum _{N\ge 0}   \sum _{  \Ga  \in
      \sL _{N,[ij]} ^k }
   \mu (\Ga )   \frac {X^N}{N!}    .
\end{equation}
and the generating series
\begin{equation}\label{Ups_glob_def}
\Upsilon (X)_{ij}=      \sum _{  \Ga  \in
      \sL _{[ij]}  }
   \mu (\Ga )   \frac {X^{N (\Ga ) }}{N(\Ga )!}    .
\end{equation}
Summing over all the trees in $\sL _{N,[ij]}^{k}$ we shall get all
the summands of the corresponding term of (\ref{chain}). Therefore
we have obtained the following formula.
\begin{prop}
\begin{equation}\label{Ups_formula}
\Upsilon _k(X)=\underbrace {HS H S H \ldots
      S H  }_{k \mbox{  times  } H}.
\end{equation}
\end{prop}
\begin{cor}
\begin{equation}\label{Ups_glob_formula}
\Upsilon (X)=H +  HS H +   HS H S H +  \ldots = H\left ( E - S H
\right ) ^{-1} = \left ( E -  H  S\right) ^{-1} H .
\end{equation}
\end{cor}

Now it is very easy to give a purely combinatorial proof of the
formula (\ref{inv_Teta}). According to (\ref{marked_half_Psi_bip})
the $ij$ component of the matrix $\Theta (S ,X)$ is the generating
function for the trees from $\bigcup _{N\ge 0} \sBB _{0,N, [ij]}$.
In any tree $\Ga \in  \sBB _{0,N, [ij]}$ there is a unique chain
connecting the two marked edges; this chain we may consider as an
element of $ \sL _{N,[ij]}^{k} $ ($k$ is the length of this chain).
Thus we have defined a mapping
\begin{equation}\label{map_c1_trees_to_chains}
   c_1 \colon  \sBB _{0,N,[ij]} \to   \sL _{N,[ij]} .
\end{equation}
As in the proof of the proposition \ref{prop_Psi_via_Psi_1} for any
chain  $\Lambda \in  \sL _{ N,[ij]}$ all the graphs in  $c_1^{-1}
(\Lambda )$ are constructed from $\Lambda $ by clutching arbitrary
genus $0$ trees with one marked tail to some of the tails of
$\Lambda $ (inserting a two-valent $s$-vertex  between the
corresponding tail of $\Ga $ and the marked edge of the tree).
Therefore
\begin{equation}\label{c1^(-1)_La}
          \sum _{\Ga  \in   c_1^{-1} (\Lambda )  }
    \frac {\mu (\Ga )}{|\Aut \Ga |}
     X^{N(\Ga )}
=     \frac {\mu (\Lambda )}{|\Aut \Lambda |}
   \left(    X +  S \Phi (S,X )\right) ^{N (\Lambda )  }
\end{equation}
and
\begin{multline}\label{proof_pro_Teta}
 \Theta ( S , X) _{ij}=   \sum _{\Ga  \in   \sBB _{0,[ij]} }
   \frac {\mu (\Ga )}{|\Aut \Ga |}   X^{N(\Ga )}
 =    \sum _{\Lambda \in   \sL _{[ij]} }
     \left( \sum _{\Ga  \in   c_1^{-1} (\Lambda )  }
    \frac {\mu (\Ga )}{|\Aut \Ga |}   X^{N(\Ga )}  \right) = \\
=  \sum _{\Lambda \in   \sL _{[ij]} }   \frac {\mu (\Lambda )}{|\Aut
\Lambda |}
   \left(    X +  S \Phi (S,X )\right) ^{N (\Lambda )  }=
   \Upsilon \left(    X +  S \Phi (S,X )\right)  _{ij}.
\end{multline}
Using (\ref{Ups_glob_formula}) and multiplying by $S$ we get
\begin{multline}\label{proof_pro_Teta1}
 S \Theta ( S , X) =S H\left(    X +  S \Phi (S,X )\right)
 +S H\left(X+S \Phi (S,X )\right)
 S H\left(X+S \Phi (S,X )\right) + \\ +
S H\left(X+S \Phi (S,X )\right)
 S H\left(X+S \Phi (S,X )\right)
S H\left(X+S \Phi (S,X )\right) + \ldots
\end{multline}
and finally adding $E$ we obtain (\ref{inv_Teta}):
\begin{equation}\label{inv_Teta_new}
E+S \Theta ( S , X) =\left ( E- S H\left(    X + S \Phi (S,X
)\right) \right) ^{-1}.
\end{equation}

This completes our study of the first term of the genus expansion
(the so-called "tree approximation"). Now let us go on with
subsequent terms. Our next step is  to describe $\Psi _g ^{(1)}$.
First let us study the case $g=1$. The set of $2$-connected genus
$1$ graphs splits into two parts $\sBB _{1,N}^{(1)} ={\sBB
'}_{1,N}^{(1)} \sqcup {\sBB ''}_{1,N}^{(1)}  $, where ${\sBB
'}_{1,N}^{(1)}$ is the set of all $2$-connected genus $1$ graphs
having only genus $0$ vertices. A connected genus $1$ graph may have
at most one vertex of positive genus; if such a vertex exists it
should have genus $1$. So a  graph $\Ga \in {\sBB ''}_{1,N}^{(1)}$
has no cycles, therefore it has no edges. Hence for each $N$ ${\sBB
''}_{1,N}^{(1)}$ consists of one graph, namely single genus $1$
vertex with $|N|$ tails colored by $N$. Therefore
\begin{equation}\label{Psi_1''}
    \sum _{N \ge 0}
     \left( \sum _{\Ga \in {\sBB ''}_{1,N}  ^{(1)}}
    \frac {\mu (\Ga )}{|\Aut \Ga |}\right)
     X^N= \Psi _1 (0,X )=U_1(X).
\end{equation}
If a  genus $1$ graph has only genus $0$ vertices then it must have
exactly one cycle. Therefore a graph $\Ga \in {\sBB '}_{1,N}^{(1)}$
consists of one cycle having $k>0$ vertices $v''\in V_a(\Ga )$
interleaving with $k$ two-valent $s$-vertices $v'\in V_s(\Ga )$.
Denote the set of such graphs by $ \sBB _{1,N}^{k\ (1)}$; the set of
such graphs with the additional choice of one  two-valent vertex
$v_0 ''\in V_s (\Ga )$ and of an orientation of the cycle will be
denoted by  $ \overline{\sBB _{1,N}^{k\ (1)}}$.
 The $2k$-sheet covering
$ \overline{\sBB _{1,N}^{k\ (1)}}\to  \sBB _{1,N}^{k\ (1)}$ is
compatible with the automorphisms of the corresponding graphs.
Therefore $\sum _{\Ga  \in  \sBB _{1,N}^{k\ (1)}   }
    \frac {\mu (\Ga )}{|\Aut \Ga  |}
     X^{N}=\frac 1{2k} \sum _{ \overline \Ga  \in
      \overline{\sBB _{1,N}^{k\ (1)}  } }
    \frac {\mu (\Ga )}{|\Aut \Ga  |}      X^N$
     where $ \overline \Ga $ denotes a graph $\Ga $ together with the
     described additional structure.
     Deletion of the vertex $v_0 ''$ defines a bijection
     \begin{equation}\label{cycle_chain}
\sBB _{1,N}^{k\ (1)} \to \bigcup _{ij} \sL _{N,[ij]}^{k} ,
\end{equation}
therefore
\begin{equation}\label{cycle_k}
\sum _{\Ga  \in  \sBB _{1,N}^{k\ (1)}   }
    \frac {\mu (\Ga )}{|\Aut \Ga  |}
     X^{N}=\frac 1{2k} \tr  \left( S\  H(X) \right) ^k.
\end{equation}

Summing for all $k$ we obtain the
   formula for $\Psi _1$    (the so-called "one-loop
   approximation").

\begin{prop}
\begin{equation}\label{Psi_1_1}
  \Psi _1  ^{(1)}(S,X )= U _1 (X) - \frac 12 \tr \ln (E-   S\ H(X)  ).
\end{equation}
\end{prop}
\begin{cor}
\begin{equation}\label{Psi_1_1_c}
    \Psi _1  (S,X )= U _1 \left( X+ S \Phi (S,X ) \right)
    - \frac 12 \tr \ln \left( E-   S H \left(  X+ S \Phi (S,X )\right) \right) .
\end{equation}
\end{cor}

A pair of two-valent $s$-vertices $v_1',v_2'    \in V_s(\Ga )$ of a
colored bipartite connected and $2$-connected  modular graph $\Ga
\in  \sBB _{g,N}^{(1)}$ with $g\ge 1$ will be called a $2$-{\it cut}
if deleting  of $v_1'$ and $v_2'$
 disconnects the graph and at least one of the two new connected components
 has genus zero. Note that for $g>1$  at most one of the two components
 may have genus zero and that the genus $0$ component is a tree
 from $\sL _{N',[ij]}^{k'}$ for some $N'$, $k'$.
A graph without $2$-cuts will be called $3$-{\it connected}; the set
of  $3$-connected genus $g$ graphs having $n_i$ tails of the color
$i$ will be denoted by
 $\sBB _{g,N}^{(2)}$;     $\sBB _{g}^{(2)}=\bigcup _{N\ge 0} \sBB _{g,N}^{(2)}$.
The $2$-cuts of a given  graph $\Ga \in  \sBB _{g,N}^{(1)}$  are
partially ordered by the inclusion relation of the corresponding
genus
 $0$ components. Consider all the
 maximal $2$-cuts.
  Replacing each corresponding maximal genus $0$ component by a new
  two-valent $s$-vertex
we obtain a $3$-connected graph $\bar \Ga \in  \sBB _{g,N}^{(2)}$.
This provides the mapping
\begin{equation}\label{map_c2_bip}
   c_2 \colon  \sBB _{g}^{(1)} \to  \sBB _{g}^{(2)}.
\end{equation}
 Pick a graph $ \Ga \in  \sBB _{g}^{ (2)} $. The
 preimage $c_2^{-1}(\bar \Ga ) $ consists
of all graphs obtained from $\Ga $ by replacing some of the
two-valent $s$-vertices by arbitrary trees from  $\sL _{N,[ij]}^{k}$
(bounded by two two-valent $s$-vertices on the clutching positions).
Therefore
\begin{equation}\label{c2^(-1)_Ga}
         \left( \sum _{\Ga  \in   c_2^{-1} (\Ga )  }
    \frac {\mu (\Ga  )}{|\Aut \Ga   |}\right)
     X^{N(\Ga   )}
\end{equation}
is obtained from $ \frac {\mu (\Ga )}{|\Aut \Ga |} X^N$ by
substituting
\begin{equation}\label{sum_SHS}
\left (  S+SHS+SHSHS+\ldots \right) _{ij}
\end{equation}
instead of all $s _{ij}$. Note that the matrix in (\ref{sum_SHS})
may be expressed as
\begin{equation}\label{sum_SHS_formula}
S+SHS+SHSHS+\ldots = S \left ( E- HS\right) ^{-1} =\left ( E-
SH\right) ^{-1}S.
\end{equation}
This enables to express the function $\Psi $ in terms of the
generating function for $3$-connected graphs
\begin{equation}\label{Psi_def_bip_2}
    \Psi _g ^{(2)}(S,X )= \sum _{N \ge 0}
     \left( \sum _{\Ga \in \sBB _{g,N}  ^{(2)}}
    \frac {\mu (\Ga )}{|\Aut \Ga |}\right)
     X^N.
\end{equation}
\begin{prop} For $g>1$
\begin{equation}\label{Psi_g_1}
    \Psi _g  ^{(1)}(S,X )=
    \Psi _g  ^{(2)}\left ( \left ( E- SH(X)\right) ^{-1} S   ,X \right).
\end{equation}
\end{prop}
\begin{cor}\label{prop_Psi_via_Psi_2} For $g>1$
\begin{equation}\label{Psi_via_Psi_2}
    \Psi _g (S,X )=
\Psi _g  ^{(2)}\left( \left ( E- SH \left(  X+ S \Phi (S,X )\right)
\right) ^{-1} S  , X+ S \Phi (S,X ) \right).
\end{equation}
\end{cor}

Now we are left to describe $\Psi _g  ^{(2)}$. Deletion of all the
tails defines the mapping
\begin{equation}\label{map_c3_bip}
   c_3 \colon  \sBB _{g}^{(2)} \to  \sBB _{g,0}^{(2)},
\end{equation}
where $ \sBB _{g,0}^{(2)}$ is the set of genus $g>1 \ $
$3$-connected graphs   without tails. Note that a graph is
$3$-connected if and only if it is stable.
 We shall  denote the
set of stable closed  graphs $\sBB _{g,0}^{(2)}$ by $\sA _g $. Pick
a stable closed graph $\Ga \in \sA _g $ and a vertex $v\in V_a(\Ga
)$; let $N(v)$ ($|N(v)|\ge 3$ for $g=0$) be the multiindex of its
valences. Note that $a_{g(v),N(v)}=\frac {\p ^{|N(v)|} U _{g(v)}}{\p
X^{N(v)}}(0)$, and the same vertex $v$ in different graphs from $c_3
^{-1} (\Ga )$ corresponds to  certain terms of the expansion of  $
\frac {\p ^{|N(v)|} U _{g(v)}(X)}{\p X^{N(v)}}$. Thus it is not hard
to verify that
\begin{equation}\label{Psi_c_3^-1}
         \left( \sum _{\De  \in   c_3^{-1} (\Ga )  }
    \frac {\mu (\De )}{|\Aut \De |}\right)
     X^{N(\De )}
=  \frac  1{|\Aut \Ga |} \prod _{v\in V_a(\Ga )}
 \frac {\p ^{|N(v)|} U _{g(v)}}{\p X^{N(v)}}(X).
\end{equation}
The product on the right side of   (\ref{Psi_c_3^-1}) looks like the
second product in the definition of $\mu (\Ga )$ (\ref{mu(Ga)}) with
$ \frac {\p ^{|N|} U _{g}(X)}{\p X^N}$ substituted instead of the
variables $a_{g,N}$. Therefore, defining for $g\ge 1$ the generating
functions $P_g$ for stable closed  graphs by
\begin{equation}\label{stable_Psi}
 P _g(\{  a_{m,N} \} ,S)=
 \sum _{\Ga  \in  \sA _g   } \frac {\mu (\Ga )}{|\Aut \Ga |}=
 \Psi _g ^{(2)}(\{  a_{m,N}  \} ,S, 0 )
\end{equation}
we obtain the following expression for $\Psi _g^{(2)}$:
\begin{multline}\label{Psi_g_2}
\Psi _g ^{(2)}\left( \left\{  a_{m,N} \right\} ,S,X \right)  = \sum
_{N \ge 0}
     \left( \sum _{\Ga \in \sBB _{g,N}  ^{(2)}}
    \frac {\mu (\Ga )}{|\Aut \Ga |}\right)
     X^N= \\
     =P_g\left( \left\{  a_{m,N}:= \frac {\p ^{|N|}  U _m }
     {\p X ^N }(X ) \right\} ,S \right) .
\end{multline}
Thus we are able to express the generating functions $\Psi _g$ in
terms of the
 generating functions for
stable closed  graphs $P_g$.
\begin{thm}\label{prop_Psi_via_P_g} For $g>1$
\begin{multline}\label{Psi_via_P_g}
    \Psi _g (S,X )= \\
=P_g\left( \left\{  a_{g,N}:=  \frac {\p ^{|N|}  U _g }
     {\p X ^N }\left( X +
S \Phi (S,X ) \right)  \right\} ,\left ( E- SH \left( X+ S \Phi (S,X
)\right) \right) ^{-1} S \right) .
\end{multline}
\end{thm}
Note that for each $g> 1$ the set of stable closed  graphs $\sA _g $
is finite: for $g>1$ a stable closed genus $g$ graph has at most
$3g-3$ two-valent $s$-vertices. Hence $ P _g(\{  a_{m,N} \} ,S )$ is
a polynomial in $ s_{ij}$ and   $\{  a_{m,N} \}$ for $ |N|\le 2g-2$
($|N|\le 2$ for $g=1$) and $0\le m\le g$. It has  degree $3g-3$ as a
polynomial in  $ s_{ij}$; for combinatorial case the degree of
 all  terms is at least $g$.  We shall call the
polynomials $P_g(\{  a_{m,N} \} ,S )$ {\it stable graph
polynomials}.
For instance the first stable graph polynomial for $g=2$ 
\begin{multline}\label{P_2}
  P_2  = a_{2,0}+\sum _{i,j} \frac 1{\{ ij \}  !}
  a_{1,\{ i \}  } a_{1,\{ j \}  } s_{ij} + \sum _{i,j} \frac 1{\{ ij \}  !}
  a_{1,\{ ij \}  }  s_{ij} + \\
  +\sum _{i,j,k,l} \frac 1{\{ kl \}  !}
  a_{1,\{ i \}  } a_{0,\{ jkl \}  } s_{ij} s_{kl}+\sum _{i,j,k,l,p,q}
  \frac 1{\Aut \Ga _5}
  a_{0,\{ ijk \}  } a_{0,\{ lpq \}  } s_{ij} s_{kl} s_{pq}+ \\
  +\sum _{i,j,k,l,p,q}
  \frac 1{\Aut \Ga _6}
  a_{0,\{ ijk \}  } a_{0,\{ lpq \}  } s_{ip} s_{kl} s_{jq}+
\sum _{i,j,k,l}  \frac 1{\Aut \Ga _7}
  a_{0,\{ ijkl \}  } s_{ij} s_{kl}.
\end{multline}
 has seven
groups of terms corresponding to the seven possible genus $2$
graphs; different terms in each group correspond to different ways
of coloring edges of the given graph. All the coefficients are the
inverses to the number of automorphisms of the corresponding graph;
in the last three terms we do not indicate an explicit expression of
dependence of  these numbers on the way of coloring. For
 combinatorial
case $P_2$ has only three last terms.

It is not hard to verify that the stable graph polynomials are
homogeneous in the following sense.
\begin{prop}\label{homogen_st_gr_poly}
Define the grading of the  polynomial ring \  $\C \left[ \{ s_{ij}\}
,\{  a_{m,N} \} \right] $ by
\begin{equation}\label{stable_deg}
\deg    a_{g,N} =   1-|N|-g \quad \mbox{and}  \quad \deg    s_{ij} =
1.
\end{equation}
Then stable graph polynomial $P_g(\{  a_{g,N}  \} ,S   ) $ is
homogeneous polynomial of degree $1-g$.
\end{prop}

\section{Stable graph polynomials}

Next let us derive the recurrence for the stable graph polynomials.
The idea of it is quite similar to the proof of theorem
\ref{th_main_1}: we delete one $s$-vertex of a given stable closed
genus $g$ graph and obtain a genus $g-1$ graph with two tails.
Unlike the cases considered in the theorems \ref{th_main_1}  the new
graph does not  correspond to the same generating function for genus
$g-1$, because it is not closed and may be not connected. First let
us study  the latter case. Let $\Ga $ be a stable closed  genus
$g>1$ graph and assume that deleting of some $s$-vertex disconnects
it into the disjoint  union of two connected graphs $\Ga '$ and $\Ga
''$. Then both $\Ga '$ and  $\Ga ''$ have positive genus and each of
the two has exactly one tail. Let us denote by $\sC _{g,\{ i \} }$
the set of all genus $g\ge 1$ graphs with the only tail of color $i$
 obtained from stable closed  graphs in the described way. Pick a graph
$\Ga \in \sC _{g,\{ i \}  }$. For $g>1$ deletion of its only tail
provides a stable closed graph unless the vertex $v_0$ incident to
the  tail was a trivalent genus $0$ vertex. But in the latter case
we obtain a stable closed  graph by substituting an $s$-vertex
instead of the subgraph consisting of the vertex $v_0$ together with
the two $s$-vertices adjacent to $v_0$ (see Fig. 2). Thus for $g>1$
we have defined a mapping
\begin{equation}\label{map_c5_bip}
   c_5 \colon  \sC _{g,\{ i \}  }  \to  \sBB _{g,0}^{(2)}.
\end{equation}
The generating function ($\mu (\Ga ) $ is defined in (\ref{mu(Ga)}))
\begin{equation}\label{stable_Q}
 Q _g ^{(i)}(\{  a_{m,N} \} ,S)=
 \sum _{\Ga  \in  \sC _{g,\{ i \}  }  } \frac {\mu (\Ga )}{|\Aut \Ga |}
\end{equation}
is a derivative of $P_g$ in the following sense. For   $1\le k\le r$
define the differentiation $D_k$ of the ring of polynomials in all
$a_{m,N}$ and $s_{ij}$ by its action on its generators:
\begin{equation}\label{der_a}
D_k( a_{m,N})=  a_{m,N+\{ k\} }
\end{equation}
\begin{equation}\label{der_s}
    D_k( s_{ij})=\sum _{p,q}  s_{ip} s_{jq} a_{0,\{ pqk \}  }
\end{equation}
\begin{prop}
The differentiation  $D_k$ is  homogeneous and has degree $-1$.
\end{prop}

Using the mapping (\ref{map_c5_bip}) it is not hard to verify the
following statement.
\begin{prop} For $g>1$
\begin{equation}\label{Q_via_P}
 Q _g ^{(i)}(\{  a_{m,N} \} ,S)= D_i( P _g (\{  a_{m,N} \} ,S)
\end{equation}
\end{prop}
For $g=1$ the mapping (\ref{map_c5_bip}) is not well-defined but it
is not hard to list all the graphs of  $\sC _{1,\{ i \}  }$
explicitly. In fact there are only two such graphs: single genus $1$
vertex with one color $i$ tail and a length one cycle with one
$s$-vertex and one genus $0$ trivalent vertex having one tail of the
color $i$. Therefore
\begin{equation}\label{Q_1}
 Q _1 ^{(i)}(\{  a_{m,N} \} ,S)=   a_{1,\{ i \}  } +
 \sum  _{p,q}  \frac 1{\{ pq \}  !}
 s _{p,q}  a_{0,\{ ipq \}  },
\end{equation}
which may be considered as a formal definition of $D_i(P_1)$
(whereas $P_1$ does not exist).

Next let us consider the second possibility. Pick a stable closed
genus $g>2$ graph and assume that deleting of some $s$-vertex does
not disconnect  it. Let us denote by $\sC _{g,\{ ij \}  }$ the set
of all genus $g> 1$ graphs having exactly two tails of colors $i$
and $j$
 obtained from stable closed  graphs in the described way.
Our next purpose is to define for $g>1$ a mapping
\begin{equation}\label{map_c4_bip}
   c_4 \colon  \sC _{g,\{ ij \}  }  \to  \sBB _{g,0}^{(2)}.
\end{equation}
First let us assume that the two tails are not attached to the same
trivalent genus $0$ vertex. In this case the mapping
(\ref{map_c4_bip}) may be described as the result of twice repeated
operations used in the definition of the mapping $c_5$
(\ref{map_c5_bip}). This corresponds to double differentiation $
\frac 1{\{ ij \}  !}D_iD_j$. Next consider a graph    $\Ga \in \sC
_{g,\{ ij \}  } $ having  two tails
  attached to the same trivalent
genus $0$ vertex  $v_0$. Then $0$ is adjacent to exactly one
$s$-vertex which connects it to the remaining part of the graph.
Removal of this $s$-vertex (together with $v_0$) provides a graph
$\bar \Ga \in \sC _{g,\{ m \}  } $ for some $m$. This enables to
describe the generating function
\begin{equation}\label{stable_R}
 R _g ^{(ij)}(\{  a_{m,N} \} ,S)=
 \sum _{\Ga  \in  \sC _{g,\{ ij \}  }  } \frac {\mu (\Ga )}{|\Aut \Ga |}
\end{equation}
using the mapping $c_4$ (\ref{map_c4_bip}).
\begin{prop} For $g>1$
\begin{equation}\label{R_via_P}
 R _g ^{(ij)}(\{  a_{m,N} \} ,S)= \frac 1{\{ ij \}  !} \left[
  D_i D_i( P _g (\{  a_{m,N} \} ,S)+
  \sum _{p,q} D_p\left( P _g (\{  a_{m,N} \} ,S)\right)
 s_{pq} a_{0,\{ ijq \}  }\right] .
\end{equation}
\end{prop}

For $g=1$ the mapping (\ref{map_c4_bip}) is not well-defined but it
is not hard to list all the graphs of  $\sC _{1,\{ ij \}  }$
explicitly. In fact there are only three such graphs corresponding
to the three terms of the following expression
\begin{multline}\label{R_1}
R _1 ^{(ij)}(\{  a_{m,N} \} ,S)= \frac 1{\{ ij \}  !} \left[ \sum
_{p,q,u,t} \frac 1{\{ \{ pq\} \{ ut\} \} ! } s _{pu}  s_{qt} a_{0,\{
ipq \} }  a_{0,\{ jut\}  }
  \right. +  \\ +   \left. \sum _{p,q} \left(
 a_{1,\{ p \}  } s_{pq} a_{0,\{ ijq\} } +
 \sum  _{u,t}  \frac 1{\{ ut \}  !}
 s _{ut}  a_{0,\{ utp \} }  s_{pq} a_{0,\{ ijq\}  } \right) \right] = \\
 = \frac 1{\{ ij \}  !} \left[
 \sum  _{p,q,u,t} \frac 1{\{ \{ pq\} \{ ut\} \} ! }
s _{pu}  s_{qt} a_{0,\{ ipq \} }  a_{0,\{ jut\}  }
 + \sum _{p,q} D_p( P _1)
 s_{pq} a_{0,\{ ijq \}  } \right] ,
\end{multline}
where $\{ \{ pq\} \{ ut\} \} !$ means $2$ for $p=q$ and $u=t$ and
$1$ otherwise. Note that in the last expression the term $\sum
_{p,q,u,t} \frac 1{\{ \{ pq\} \{ ut\} \} ! } s _{pu}  s_{qt} a_{0,\{
ipq \} }  a_{0,\{ jut\}  } $ may be considered as a formal
definition of $D_iD_j(P_1)$ (different from  $D_i(D_j(P_1))$  and
$D_j(D_i(P_1))$ which are not equal).

Now we are prepared  to present the recurrence for stable graph
polynomials.
 Deletion of an  $s$-vertex from a given stable closed
genus $g$ graph provides either a connected graph from $\sC _{g-1,\{
ij \}  }$ or a pair of connected graph from $\sC _{m,\{ i \}  }$ and
$\sC _{g-m,\{ j \}  }$ for some $1\le m\le g-1$. Therefore
\begin{equation}\label{rec_predv}
    \frac {\p Pg}{\p s_{ij}}= \frac 1{\{ ij \}  !} \left[
R _{g-1} ^{(ij)}(\{  a_{m,N} \} ,S)+
  \sum _{m=1}^{g-1} Q_m ^{(i)}(\{  a_{m,N} \} ,S)
 Q_{g-m} ^{(j)}(\{  a_{m,N} \} ,S)
   \right]  .
\end{equation}
Substituting (\ref{Q_via_P}),   (\ref{Q_1}) and   (\ref{R_via_P}) we
get the desired recurrence.
\begin{thm}\label{prop_rec_P_g}   For $g>2$ stable graph
polynomials $P_g(\{  a_{m,N} \} ,S)$  for colored  bipartite graphs
are given by the recurrences (for each pair $ij$)
\begin{multline}\label{rec_P_g}
    \frac {\p Pg}{\p s_{ij}}= \frac 1{\{ ij \}  !} \left[
  D_i D_j( P _{g-1})+
  \sum _{p,q} D_p( P _{g-1})
 s_{pq} a_{0,\{ ijq \}  }+ \right. \\
+ D_i( P _{g-1})\left( a_{1,\{ j \}  } + \sum  _{p,q}  \frac 1{\{ pq
\}  !}
 s _{p,q}  a_{0,\{ jpq \}  }
\right) + \\
+\left. D_j( P _{g-1})\left( a_{1,\{ i \}  } + \sum  _{p,q}  \frac
1{\{ pq \}  !}
 s _{p,q}  a_{0,\{ ipq \}  }
\right) +
 \sum _{m=2}^{g-2} D_i ( P _{m})D_j( P _{g-m})
 \right]  .
\end{multline}
and by the initial condition $P_g(\{  a_{m,N} \} ,0)= a_{g,0}$,
where the differentiations $D_i$ are given by (\ref{der_a}) and
(\ref{der_s}) and $P_2$ is defined in  (\ref{P_2}) .
\end{thm}

Regardless of the absence of genus $1$ stable closed  graphs we may
start the recurrence (\ref{rec_P_g}) from $g=1$ formally putting
\begin{equation}\label{D_i(P_1)}
D_i(P_1)= a_{1,\{ i \}  } + \sum  _{p,q}  \frac 1{\{ pq \}  !}
 s _{p,q}  a_{0,\{ ipq \}  } .
\end{equation}
Then for any $g>1$
\begin{multline}\label{rec_P_g_1}
    \frac {\p Pg}{\p s_{ij}}= \frac 1{\{ ij \}  !} \left[
  D_i D_j( P _{g-1})+
  \sum _{p,q} D_p( P _{g-1})
 s_{pq} a_{0,\{ ijq \}  }+ \right. \\
 +\left.
 \sum _{m=1}^{g-1} D_i ( P _{m})D_j( P _{g-m})
 \right] .
\end{multline}


\section{Counting functions for combinatorial graphs.} \label{r=1}

For $r=1$ we have one variable $s$ and two-index variables
$a_{g,n}$. The differentiation $D$ of the ring  of polynomials in
$s$ and $a_{g,n}$ is defined by
\begin{equation}\label{der_a_1}
D( a_{g,n})=  a_{g,n+1 }
\end{equation}
\begin{equation}\label{der_s_1}
    D( s)=  s^2 a_{0,3 }
\end{equation}
The polynomial $P_2$ is given by (\ref{P_2}):
\begin{multline}\label{P_2_1}
  P_2  = a_{2,0}+ \frac 12
  a_{1,1 } ^2  s +  \frac 12
  a_{1,2 }  s
  +\frac 12
  a_{1,1  } a_{0,3  } s^2+\frac 18 a_{0,4  } s^2+
 \left( \frac 1{12} + \frac 18 \right)
  a_{0,3  } ^2 s^3,
\end{multline}
and for combinatorial case ($a_{g,n}=0$ for $g>0$)
\begin{equation}\label{P_2_1_comb}
  P_2  = \frac 18 a_{0,4  } s^2+
  \frac 5{24}
  a_{0,3  } ^2 s^3.
\end{equation}
The recurrence (\ref{rec_P_g}) becomes
\begin{multline}\label{rec_P_g_r=1}
    \frac {\p Pg}{\p s}= \frac 12 \left[
  D^2( P _{g-1})
   +2 D( P _{g-1})\left( a_{1,1  } +
 s   a_{0,3 }
\right) +
 \sum _{m=2}^{g-2} D( P _{m})D( P _{g-m})
 \right]
\end{multline}
and the formula (\ref{Psi_via_P_g}) looks like

\begin{multline}\label{Psi_via_P_g_r=1}
    \Psi _g (s,x )= \\
=P_g\left( \left\{   a _{g,n}: =  \frac {d ^n  U _g }
     {d x ^n }\left( x +
s \Phi (s,x ) \right)  \right\} ,\frac s{1- s  H \left(  x+ s \Phi
(s,x )\right)}
      \right) .
\end{multline}


For the  counting functions for   combinatorial graphs  we put
\begin{equation}\label{all_gr_mu}
    a _{g,n}^{\mathrm{comb}} = \left\{\begin{array}{cc}
      1  & \mbox{ if  }g=0  \\
      0 & \mbox{ otherwise }  \\
     \end{array}\right. \end{equation}
which we may write uniformly as $  a _{g,n}^{\mathrm{comb}}  = \de
_{g0}$. Then the initial conditions (see (\ref{Psi_init})
  --- (\ref{F_def})) are
\begin{equation}\label{all_gr_F}\begin{array}{c}
   U_g^{\mathrm{comb}}(x)=0  \quad \mbox{    for    }  \quad g>0, \\
  U_0^{\mathrm{comb}}(x)=e^x, \\
  \end{array}
 \end{equation}
 and hence $F^{\mathrm{comb}}(x)=H^{\mathrm{comb}}(x)=e^x$.
The counting series for all trees satisfies the functional equation
(\ref{fun_eq_Phi_bip}):
\begin{equation}\label{all_gr_fun_eq}
\Phi ^{\mathrm{comb}}(s,x)=e^{ x+ s \Phi ^{\mathrm{comb}}(s,x)}
 \end{equation}
 and therefore
\begin{equation}\label{all_gr_U_g}
    \frac {d ^n  U _g^{\mathrm{comb}} }
     {d x ^n }\left( x +
s \Phi ^{\mathrm{comb}}(s,x ) \right)=0 \quad \mbox{    for    }
\quad g>0,
 \end{equation}
and
\begin{equation}\label{all_gr_U_0}
    \frac {d ^n  U _0^{\mathrm{comb}} }
     {d x ^n }\left( x +
s \Phi ^{\mathrm{comb}}(s,x ) \right)=e^{ x+ s \Phi
^{\mathrm{comb}}(s,x)} =  \Phi ^{\mathrm{comb}}(s,x) \quad \mbox{
for    } \quad g=0,
 \end{equation}
which we may write in a  uniform way  as $  a _{g,n}: =  \Phi
^{\mathrm{comb}}(s,x)\de _{g0}$. Denote  the second argument of
(\ref{Psi_via_P_g_r=1}) by $Y$:
\begin{equation}\label{all_gr_sH}
 Y=   \frac {s}
     {1-s H\left( x +
s \Phi (s,x ) \right)}:= \frac {s}
     {1-s e^{ x +
s \Phi ^{\mathrm{comb}}(s,x ) }} = \frac {s}
     {1-s \Phi ^{\mathrm{comb}}(s,x ) }
.
 \end{equation}

Therefore for the counting function for combinatorial graphs
\begin{multline}\label{Psi_via_P_g_comb}
    \Psi _g^{\mathrm{comb}} (s,x )= \Psi _g (\left\{   a _{g,n}: =
    \de _{g0}  \right\} ,  s,x )= \\
=P_g\left( \left\{   a _{g,n}: =
 \Phi ^{\mathrm{comb}}(s,x )\de _{g0}  \right\} ,
\frac s{1- s \Phi ^{\mathrm{comb}}(s,x ) }
 \right) .
\end{multline}

Recall that stable graph polynomials are homogenous and have degree
$1-g$ (see proposition \ref{homogen_st_gr_poly}) with respect to the
grading (\ref{stable_deg}). By definition of stable graph
polynomials (\ref{stable_Psi}):
\begin{equation}\label{stable_P_g_r=1}
 P _g(\{  a_{m,N} \} ,s)=
 \sum _{\Ga  \in  \sA _g } \frac {\mu (\Ga )}{|\Aut \Ga |} ,
\end{equation}
where $\sA _g ^k$ is the set of genus $g$ stable closed  graphs.
Denote the set of genus $g$  stable closed  graphs with $k$ edges by
$\sA _g ^k$.
 A graph $\Ga \in \sA _g ^k$ has $k-g+1$ vertices (combinatorial
 graphs have no vertices of higher genus), so  $\mu (\Ga ) = s^k
\prod _{i=1}^{k-g+1} a_{0,n_i}$, where $n_i$ are the valences of the
vertices. Therefore the counting function for combinatorial
 graphs
\begin{multline}\label{stable_P_g_r=1_1}
 \Psi _g ^{\mathrm{comb}}(s,x )
= P _g\left( \left\{  a_{m,N}: =
 \Phi ^{\mathrm{comb}}(s,x )\de _{m0}  \right\} ,Y\right) =\\
 = \sum _k Y^k
\left( \sum _{\Ga  \in  \sA _g ^k} \frac 1{|\Aut \Ga |}\prod
_{i=1}^{k-g+1} a_{0,n_i}\right)
=  \sum _k Y^k \left( \sum _{\Ga  \in  \sA _g ^k}
\frac 1{|\Aut \Ga |} \Phi ^{\mathrm{comb}}(s,x )^{k-g+1}\right) = \\
=\frac 1{\Phi (s,x )^{g-1}} \sum _k  \left(  \Phi
^{\mathrm{comb}}(s,x )Y\right)  ^k \left( \sum _{\Ga  \in  \sA _g
^k}
\frac 1{|\Aut \Ga |}\right) = \\
 =\frac 1{\Phi ^{\mathrm{comb}}(s,x )^{g-1}}  P _g\left( \left\{  a_{m,N}:
 =  \de _{m0}  \right\} ,
 \frac {s\Phi ^{\mathrm{comb}}(s,x )}{1- s \Phi ^{\mathrm{comb}}(s,x ) }\right).
\end{multline}
Note that the polynomials
\begin{equation}\label{P_g^comb}
   P_g^{\mathrm{comb}}(s)= P _g\left( \left\{  a_{m,N}: =  \de _{m0}  \right\} ,
 s \right)
\end{equation}
are the generating functions for combinatorial stable closed
graphs. Thus we have proved the second part of the theorem
\ref{count_comb_st}.

Next let us prove the first part of this theorem.
 The only problem in deriving an explicit recurrence for $P_g^{\mathrm{comb}}(s)$ from
  (\ref{rec_P_g})   is
how to express the result of substitution $  a_{m,N}: =  \de _{m0}$
into $D(P_g)$ and $D^2(P_g)$ in terms of $P_g^{\mathrm{comb}}(s)$.

Consider the polynomial ring $\C [\al ,Y]$ with the grading
\begin{equation}\label{deg_Y_al}
\deg Y=1 \quad \mbox{ and}  \quad \deg \al  =-1
\end{equation}
 and differentiation $\de ^{\mathrm{comb}}$ of this ring
defined by its action on the generators
\begin{equation}\label{diff_Y_al}
\de _{\mathrm{comb}}(\al )=\al \quad \mbox{ and}  \quad
 \de _{\mathrm{comb}}(Y )=Y^2 \al .
\end{equation}
Then $\de _{\mathrm{comb}}$ is homogeneous of degree $0$. Define a
ring homomorphism
\begin{equation}\label{f_comb}
f ^{\mathrm{comb}} \colon \C \left[ \{ a_{g,n} \},s_{ij} ) \right]
\to \C [\al ,Y]
\end{equation}
 by its action on the generators:
\begin{equation}\label{f_comb_a_s}
\begin{array}{l}
   f ^{\mathrm{comb}} ( a_{g,n} )=0  \quad \mbox{ for}  \quad  g>0, \\
  f ^{\mathrm{comb}} ( a_{0,n} )=\al \\
  f ^{\mathrm{comb}} (s_{ij}) =Y. \\
 \end{array}
\end{equation}
Evidently $D \circ f ^{\mathrm{comb}} = f ^{\mathrm{comb}} \circ
\de _{\mathrm{comb}}$, therefore
 $D^2 \circ f ^{\mathrm{comb}}
= f ^{\mathrm{comb}} \circ  \de _{\mathrm{comb}}^2$. The calculation
(\ref{stable_P_g_r=1_1}) shows that $ f ^{\mathrm{comb}}(P_g)$ is
homogeneous polynomial of degree $g-1$, hence the same is true about
$ f ^{\mathrm{comb}}(D(P_g)))= \de _{\mathrm{comb}} ( f
^{\mathrm{comb}}(P_g))$. But for any degree $g-1$ homogeneous
polynomial $W\in \C [\al ,Y]$
\begin{multline}\label{de_comb}
\de _{\mathrm{comb}}(W) =Y^2 \al \frac {\p W}{\p Y }+ \al \frac {\p
W}{\p \al } = \\ =Y^2 \al \frac {\p W}{\p Y }+Y \frac {\p W}{\p Y }-
Y \frac {\p W}{\p Y }+\al \frac {\p W}{\p \al } = \\
=Y (Y \al +1)\frac {\p W}{\p Y }- (g-1)W.
\end{multline}
(On the last step we  use the Euler formula for homogeneous
polynomials: $Y \frac {\p W}{\p Y }-\al \frac {\p W}{\p \al
}=(g-1)W.$) Now we can apply (\ref{de_comb}) to the combinatorial
stable graph polynomials
\begin{equation}\label{P_comb=f(P)}
    P_g^{\mathrm{comb}}(s) = (f ^{\mathrm{comb}} P_g)(1,s)
\end{equation}
and get the following recurrence for $P_g^{\mathrm{comb}}$.
\begin{prop}
The combinatorial stable graph polynomials  $P_g^{comb}$ (see
(\ref{P_g^comb}))
 satisfy
the following recurrence
\begin{multline}\label{rec_P_g_r=1_count}
    \frac {d P_g^{\mathrm{comb}}}{d s}= \frac 12 \left[
  D_{comb}^2( P _{g-1}^{\mathrm{comb}})
   +2s D_{comb}( P _{g-1}^{\mathrm{comb}}) +
 \sum _{m=2}^{g-2} D_{comb}( P^{\mathrm{comb}} _{m})D_{comb}( P^{\mathrm{comb}} _{g-m})
 \right] ,
\end{multline}
where
\begin{equation}\label{D_comb}
   D_{comb}= s(s+1)\frac{d}{ds}  - (g-1).
\end{equation}
\end{prop}

Here we present the explicit form of the recurrence
\ref{rec_P_g_r=1_count} and  the polynomials $P_g$ for $g\le 6$:
\begin{multline}\label{rec_P_g_r=1_count_expl}
    \frac {d }{d s} P_g^{\mathrm{comb}}(s)= \frac 12 \left[
s^2(s+1)^2\frac{d^2}{ds^2} P_{g-1}^{\mathrm{comb}}(s)-\right.    \\
-
s(s+1)(2g-4s-5)\frac{d}{ds} P_{g-1}^{\mathrm{comb}}(s) + (g-1) P_{g-1}^{\mathrm{comb}}(s) +\\
 +
 \sum _{m=2}^{g-2} \left(
\left( s(s+1)\frac{d}{ds} P_m^{\mathrm{comb}}(s) -
(m-1)P_m^{\mathrm{comb}}(s)\right) \times   \right.  \\  \left.
\left.  \times \left( s(s+1)\frac{d}{ds} P_{g-m}^{\mathrm{comb}}(s)
- (g-m-1)P_{g-m}^{\mathrm{comb}}(s)\right)
 D( P _{m})D( P _{g-m})
 \right) \right]
\end{multline}
Here we present the polynomials $P_g^{\mathrm{comb}}$ for $g\le 6$:
$$P _2^{\mathrm{comb}}={\frac {5}{24}}\,s^3+\frac 18 s^2$$
$$P_3^{\mathrm{comb}}
={\frac {5}{16}}\,{s}^{6}+{\frac {25}{48}}\,{s}^{5}+{\frac
{11}{48}}\,s^{4} +\frac 1{48}s^{3}
$$
$$P_4^{\mathrm{comb}}
={\frac {1105}{1152}}\,{s}^{9}+{\frac {985}{384}}\,{s}^{8}+{\frac
{1373 }{576}}\,{s}^{7}+{\frac {515}{576}}\,{s}^{6}+{\frac
{223}{1920}}\,s^{5}+{ \frac {1}{384}}s^{4}
$$
$$P_5^{\mathrm{comb}}
={\frac {565}{128}}\,{s}^{12}+{\frac {12455}{768}}\,{s}^{11}+{\frac
{ 26581}{1152}}\,{s}^{10}+{\frac {12227}{768}}\,{s}^{9}+{\frac
{2089}{384 }}\,{s}^{8}+{\frac {9583}{11520}}\,{s}^{7}+{\frac
{27}{640}}\,s^{6}+{ \frac {1}{3840}}s^{5}
$$
$$P_6^{\mathrm{comb}}
={\frac {82825}{3072}}\,{s}^{15}+{\frac
{387005}{3072}}\,{s}^{14}+{\frac {371195}{1536}}\,{s}^{13}+{\frac
{10154003}{41472}}\,{s}^{12}+{\frac { 121207}{864}}\,{s}^{11}+
$$
$$+{\frac {519883}{11520}}\,{s}^{10}+{\frac {
1573507}{207360}}\,{s}^{9}+{\frac {2597}{4608}}\,{s}^{8}+{\frac
{803}{ 64512}}\,s^{7}+{\frac {1}{46080}}s^{6}$$

\bigskip

Similar formulas describe counting functions for combinatorial
stable graphs. For this case
\begin{equation}\label{all_st_gr_mu}
    a _{g,n}^{\mathrm{st}} = \left\{\begin{array}{cc}
      1  & \ \mbox{ if  } \ g=0   \  \mbox{ and  } n\ge 3 \\
      0 & \mbox{ otherwise, }  \\
     \end{array}\right. \end{equation}
so the initial conditions 
  are
\begin{equation}\label{all_st_gr_F}\begin{array}{c}
   U_g^{\mathrm{st}}(x)=0  \quad \mbox{    for    }  \quad g>0, \\
  U_0^{\mathrm{st}}(x)=e^x-1-x-\frac{x^2}2, \\
  \end{array}
 \end{equation}
 and hence $F^{\mathrm{st}}(x)=e^x-1-x$ and   $H^{\mathrm{st}}(x)=e^x-1$.
The counting series for all stable trees satisfies the functional
equation (\ref{fun_eq_Phi_bip}):
\begin{equation}\label{all_gr_st_fun_eq}
\Phi ^{\mathrm{st}}(s,x)=e^{ x+ s \Phi ^{\mathrm{st}}(s,x)}-1- (x+ s
\Phi ^{\mathrm{st}}(s,x))
 \end{equation}
 and therefore
\begin{equation}\label{all_gr_U_g}
    \frac {d ^n  U _g^{\mathrm{st}} }
     {d x ^n }\left( x +
s \Phi ^{\mathrm{st}}(s,x ) \right)=0 \quad \mbox{    for    } \quad
g>0,
 \end{equation}
and
\begin{equation}\label{all_gr_U_0}
    \frac {d ^n  U _0^{\mathrm{st}} }
     {d x ^n }\left( x +
s \Phi ^{\mathrm{st}}(s,x ) \right)=e^{ x+ s \Phi
^{\mathrm{st}}(s,x)} = 1+x+(s+1) \Phi ^{\mathrm{st}}(s,x) \quad
\mbox{ for } \quad g=0, \ n\ge 3.
 \end{equation}
 Since the terms $\frac {d ^n  U _0^{\mathrm{st}} }
     {d x ^n }\left( x +
s \Phi ^{\mathrm{st}}(s,x ) \right)$ for $n<3$ are not involved in
the stable graph polynomials  for the use in formula
(\ref{Psi_via_P_g_r=1})
 we may write in  a  uniform way
$  a _{g,n}: = ( 1+x+(s+1) \Phi ^{\mathrm{st}}(s,x) )\de _{g0}$.

The second argument of (\ref{Psi_via_P_g_r=1}) for this case is :
\begin{equation}\label{all_gr_sH}
   \frac {s}
     {1-s H^{\mathrm{st}}\left( x +
s \Phi ^{\mathrm{st}}(s,x ) \right)}= \frac {s}
     {1-s (e^{ x +
s \Phi ^{\mathrm{st}}(s,x ) }-1)} = \frac {s}
     {1-s (x+(s+1)\Phi ^{\mathrm{st}}(s,x ) )}
.
 \end{equation}
Thus using the same argument as in the proof of
(\ref{stable_P_g_r=1_1}) we get the formula for counting function
for combinatorial stable graphs:
\begin{multline}\label{stable_P_g_r=1_st}
 \Psi _g ^{\mathrm{st}}(s,x )= \\
= P _g\left( \left\{  a_{m,N}: =
 \left( 1+x+(s+1) \Phi ^{\mathrm{st}}(s,x) \right)  \de _{m0}  \right\} ,\frac {s}
     {1-s (x+(s+1)\Phi ^{\mathrm{st}}(s,x ) )}\right) =\\
 =\frac 1{\left( 1+x+(s+1) \Phi ^{\mathrm{st}}(s,x) \right)^{g-1}}
  P _g\left( \left\{  a_{m,N}:
 =  \de _{m0}  \right\} ,
 \frac {s\left( 1+x+(s+1) \Phi ^{\mathrm{st}}(s,x) \right)}
 {1-s (x+(s+1)\Phi ^{\mathrm{st}}(s,x ) ) }\right)   = \\
 =\frac 1{\left( 1+x+(s+1) \Phi ^{\mathrm{st}}(s,x) \right)^{g-1}}
  P _g  ^{\mathrm{comb}}  \left(
 \frac {s\left( 1+x+(s+1) \Phi ^{\mathrm{st}}(s,x) \right)}
 {1-s (x+(s+1)\Phi ^{\mathrm{st}}(s,x ) ) }\right)
 .
\end{multline}

The third  part of the theorem  \ref{count_comb_st} is proved.

\end{document}